\documentclass[notitlepage]{jams-l} 
\usepackage{amsthm,amsmath,amssymb,latexsym,epsfig} 

\newlength{\jmr}
\newlength{\hwl}
\newlength{\khov}
\newlength{\bernd}
\newlength{\fw}
\newtheorem{rolle}{Rolle's Theorem} 

\newtheorem{kho}{Khovanski's Theorem on Real Fewnomials (Special Case)} 

\newtheorem{dime}{Real Dimension Lemma} 

\newtheorem{nota}{Notation}

\newtheorem{lemma}{Lemma}
\newtheorem{prop}{Proposition}
\newtheorem{dfn}{Definition}

\newtheorem{main}{theorem} 
\newtheorem{thm}[main]{Theorem}
\newtheorem{gkconj}{Generalized Kushnirenko's Conjecture (GKC)}

\newtheorem{des}{Univariate Generalized Descartes' Rule of Signs (UGDRS)}

\newtheorem{cor}{Corollary}
\newtheorem{rem}{Remark}

\newcommand{\eps}{\varepsilon}

\newcommand{\cK}{\mathcal{K}}
\newcommand{\cN}{\mathcal{N}}

\newcommand{\supp}{\mathrm{Supp}}

\newcommand{\conv}{\mathrm{Conv}}
\newcommand{\newt}{\mathrm{Newt}}

\newcommand{\thth}{{\underline{\mathrm{th}}}}
\newcommand{\rd}{ {\underline{ \mathrm{rd} } }  }

\newcommand{\nd}{{\underline{\mathrm{nd}}}}

\newcommand{\R}{\mathbb{R}}
\newcommand{\C}{\mathbb{C}}
\newcommand{\N}{\mathbb{N}}
\newcommand{\Z}{\mathbb{Z}}

\newcommand{\Rn}{\R^n}

\newcommand{\Cn}{\C^n}

\newcommand{\Rs}{\R^*}

\renewcommand{\qed}{$\blacksquare$}
\newcommand{\cL}{{\mathcal{L}}}

\newcommand{\bO}{\mathbf{O}}

\begin{document}

\title[Trinomials in the Plane and Beyond]{Counting Isolated Roots of 
Trinomial Systems in the Plane and Beyond}


\author{Tien-Yien Li}
\author{J.\ Maurice Rojas}\thanks{
Li was partially supported by a Guggenheim Fellowship. 
Rojas' work on this paper was partially supported by Hong Kong UGC Grant 
\#9040402-730, Hong Kong/France PROCORE Grant \#9050140-730, 
and a grant from the Texas A\&M University Faculty of Science.
Some of Wang's research was done during a stay at the
Key Laboratory for Symbolic Computation and Knowledge
Engineering of the Ministry of Education, P.\ R.\ China. Wang's research 
is supported in part by the Visiting Scholar
Foundation of Key Labs In Universities, Ministry of
Education, P.\ R.\ China. } 
\author{Xiaoshen Wang}
\address{Department of Mathematics, Michigan State University,
            East Lansing, Michigan \ 48824, USA.}
\email{li@math.msu.edu } 
\urladdr{http://www.mth.msu.edu/\~{}li } 

\address{Department of Mathematics, Texas A\&M University, College Station, 
Texas \ 77843-3368, USA. }  
\email{rojas@math.tamu.edu } 
\urladdr{http://www.math.tamu.edu/\~{}rojas } 

\address{Department of Mathematics and Statistics, University 
of Arkansas at Little Rock, Little Rock, Arkansas \ 72204, USA.}  
\email{xxwang@ualr.edu } 
\urladdr{http://dms.ualr.edu/Faculty/Wang.html } 

\date{\today} 

\mbox{}\\
\vspace{-.5in} 
\maketitle

\mbox{}\\
\vspace{-1.5cm} 
\begin{abstract} 
We prove that any pair of bivariate trinomials has at most 
$5$ isolated roots in the positive quadrant. The best previous upper 
bounds independent of the polynomial degrees counted only non-degenerate 
roots and even then gave much larger bounds, e.g., $248832$ via 
a famous general result of Khovanski. Our bound is sharp, allows real 
exponents, and extends to certain systems of $n$-variate fewnomials, 
giving improvements over earlier bounds by a factor exponential in the 
number of monomials. We also derive new sharper bounds on the number of 
real connected components of fewnomial hypersurfaces. 
\end{abstract} 

\section{Introduction} 
Generalizing Descartes' Rule of Signs to polynomial systems has 
proven to be a significant challenge. Recall that a weak version of this 
famous classical result asserts that any real univariate polynomial with 
exactly $m$ monomial terms has at most $m-1$ positive roots. This bound is 
sharp and generalizes easily to real exponents (cf.\ section \ref{sec:back}). 
The original statement in Ren\'e Descartes' {\it La G\'eom\'etrie}\/ pre-dates 
1641.  Proofs can be traced back to work of Gauss in 
1828 and other authors earlier, but a definitive sharp bound for multivariate 
polynomial systems seems to have elluded us in the second millenium. This is 
particularly unfortunate since sparse polynomial systems now 
occur in applications as diverse as radar imaging 
\cite{keith} and chemistry \cite{gaterhub}. 

One simple way to generalize the setting of Descartes' Rule to higher 
dimensions and {\bf real} exponents is the following: 
\begin{nota} 
For any $c\!\in\!\Rs\!:=\R\!\setminus\!\{0\}$ and 
$a\!=\!(a_1,\ldots,a_n)\!\in\!\Rn$, let 
$x^a\!:=\!x^{a_1}_1\cdots x^{a_n}_n$ and call $c x^a$ a 
{\bf monomial term}. We will refer to $\R^n_+\!:=\!\{x\!\in\!\Rn \; | 
\; x_i\!>\!0 \text{ \ for \ all \ } i \}$ as the positive {\bf orthant}, or 
{\bf quadrant} or {\bf octant} when $n$ is respectively $2$ or $3$. 
Henceforth, we will 
assume that $\pmb{F}\!:=\!(f_1,\ldots,f_k)$ where, for all $i$, 
$f_i\!\in\!\R[x^a \; | \; a\!\in\!\Rn]$ and $f_i$ has exactly $m_i$ monomial 
terms. We call $f_i$ an {\bf $\pmb{n}$-variate $\pmb{m_i}$-nomial}\footnote{
Quite naturally, we will also call $2$-nomials {\bf binomials} and $3$-nomials 
{\bf trinomials}.} and, when $m_1,\ldots,m_k\!\geq\!1$, we call $F$ a 
{\bf $\pmb{k\times n}$ fewnomial system\footnote{We use this terminology 
solely for succinctness. 
{\bf Fewnomial theory} \cite{few} is an important related body of work 
regarding a special class of functions which includes our $m$-nomials.}   
 (over $\pmb{\R}$) of type} $\pmb{(m_1,\ldots,m_k)}$.  
Finally, we say a real root $\zeta$ of $F$ is {\bf isolated} 
(resp.\ {\bf non-degenerate}) iff the only arc\footnote{i.e., 
point or homeomorphic image of the unit circle or (open, closed, or half-open) 
unit interval...} of real roots of $F$ 
containing $\zeta$ is $\zeta$ itself \mbox{(resp.\ the Jacobian of $F$, 
evaluated at $\zeta$, has full rank). \qed} 
\end{nota} 
\begin{gkconj} 
Suppose $F$ is an $n\times n$ 
fewnomial system of type $(m_1,\ldots,m_n)$. Then the maximum number of 
non-degenerate roots of $F$ in the positive orthant is $\prod^n_{i=1} 
(m_i-1)$. \qed 
\end{gkconj} 
\begin{rem} 
\label{rem:pyr} 
The polynomial system 
$\left(\prod^{m_1-1}_{i=1}(x_1-i),\ldots, \prod^{m_n-1}_{i=1}(x_n-i)\right)$  
easily shows that the conjectured 
maximum can at least be attained (if not exceeded), and integral exponents and 
coefficients suffice for this to happen. $\diamond$ 
\end{rem} 
We can then succinctly state the original Kushnirenko's Conjecture (formulated 
in the mid-1970's by Anatoly G.\ Kushnirenko) as the 
special case of GKC where all the exponents of $F$ are non-negative integers. 
Curiously, Kushnirenko's Conjecture was open for nearly three decades until 
Bertrand Haas found a counter-example in the case 
$(n,m_1,m_2)\!=\!(2,3,3)$ (see remark \ref{rem:haas} below and \cite{haas}). 
So we will derive a correct and sharp extension of Descartes' bound to this 
case, as well as certain additional cases with $n\!\geq\!3$, $m_n\!\geq\!4$, 
and degeneracies allowed. Interestingly, the introduction of real exponents 
and degeneracies gives us more flexibility than trouble: The 
proof of our first main result uses little more than exponential 
coordinates and Rolle's Theorem from calculus. 
\begin{dfn}
For any $m_1,\ldots,m_n\!\in\!\N$, let $\pmb{\cN(m_1,\ldots,m_n)}$ 
(resp.\ $\pmb{\cN'(m_1,\ldots,m_n)}$) denote the maximal number of isolated 
(resp.\ non-degenerate) roots an $n\times n$ fewnomial system of type 
$(m_1,\ldots,m_n)$ can have in the positive orthant. \qed 
\end{dfn} 
\begin{thm} 
\label{thm:tri3} 
For all $m\!\geq\!3$ we have 
$\cN(3,m)\!\leq\!2^m-2$ and, in particular, $\cN(3,3)\!=\!5$, 
$\cN(3,4)\!\leq\!14$, and $\cN(3,5)\!\leq\!30$. 
Furthermore, $\cN'(3,m)\!=\!\cN(3,m)$. 
\end{thm} 
\noindent 
The quantities $\cN'(1,m_2,\ldots,m_n)$, $\cN(1,m_2,\ldots,m_n)$, 
$\cN'(2,m_2)$, and $\cN(2,m_2)$ are much easier to to compute than $\cN(3,3)$:  
explicit formulae for them are stated in theorem \ref{thm:tri1} of section 
\ref{sec:back}. 
\begin{rem} 
The value of $\cN(3,3)$ was previously unknown and the authors 
are unaware of any earlier result implying the equality 
$\cN'(3,m)\!=\!\cN(3,m)$. In particular, the only other information 
previously known about $\cN'(3,m)$ or $\cN(3,m)$ was an upper bound of 
$3^{m+2}2^{(m+2)(m+1)/2}$ for $\cN'(3,m)$ (see remark \ref{rem:kho} below). 
For $m\!=\!3,4,5$ the latter formula evaluates to $248832$, $23887872$, and 
$4586471424$ respectively. $\diamond$ 
\end{rem} 
\begin{rem}
\label{rem:degen} 
Note that $\cN'(m_1,\ldots,m_n)\!\leq\!\cN(m_1,\ldots,m_n)$ for all 
$m_1,\ldots,m_n\!\in\!\N$, since non-degenerate roots of $n\times n$ 
fewnomial systems are always isolated roots. While we do not yet know 
of any cases where the inequality is strict, it is interesting to note 
that GKC can {\bf not} be strengthened to allow degeneracies: For example, 
the polynomial system\footnote{Examples of this type were observed earlier 
by William Fulton around 1984 (see the first edition of \cite
{ifulton}) and Bernd Sturmfels around 1997 \cite{poly}.}
\scalebox{.9}[1]{$\left(x_1(x_3-1),x_2(x_3-1), 
\prod^5_{i=1}(x_1-i)^2+\prod^5_{i=1}(x_2-i)^2 \right)$}  
is of type $(2,2,21)$, has $25$ integral roots in the positive octant (all of 
which have singular Jacobian), but its GKC bound is $20$.  $\diamond$ 
\end{rem}
\begin{rem} 
\label{rem:haas} 
Haas' counter-example to the original Kushnirenko's Conjecture is 
\[(x^{108}_1+1.1x^{54}_2-1.1x_2, x^{108}_2+1.1x^{54}_1-1.1x_1),\]  
which has $5$ roots in the positive quadrant, thus contradicting its alleged 
GKC bound of $4$ \cite{haas}.\footnote{Dima Grigoriev informed the 
author on Sept.\ 8, 2000 that Konstantin A.\ Sevast'yanov, a colleague of 
Kushnirenko and contemporary of Grigoriev, had found a similar counter-example 
much earlier. Unfortunately, this counter-example does not seem to have been 
recorded and, tragically, Sevast'yanov committed suicide some time \mbox{before 
1997.}}  Jan Verschelde has also verified numerically that there are 
exactly $108^2\!=\!11664$ complex roots, and thus (assuming the floating-point  
calculations were sufficiently good) each root is non-degenerate by 
\mbox{B\'ezout's theorem. $\diamond$ }  
\end{rem}

The central observation that led to our proof may be of independent interest. 
We state it as assertion (3) of theorem \ref{thm:cool} below. The first 
two assertions dramatically refine the bounds of 
Oleinik, Petrovsky, Milnor, Thom, and Basu on the number of connected 
components of a real algebraic set \cite{op,milnor,thom,basu} in the 
special case of a single polynomial and extend to real exponents: 
\begin{thm} 
\label{thm:cool} 
Let $\pmb{Z}$ be the set of roots in $\R^n_+$ of an $n$-variate $m$-nomial. 
Also let $\pmb{\cK'(n,\mu)}$ denote the maximal number of non-degenerate 
roots in $\Rn_+$ of an $n\times n$ fewnomial system with exactly $\mu$ 
distinct exponent vectors. Finally, let 
$P_{\mathrm{comp}}(n,m)$ (resp.\ $P_{\mathrm{non}}(n,m)$)  
be the maximal number of compact (resp.\ non-compact) connected components of 
any such $Z$. Then... 
\begin{enumerate} 
\item{$P_{\mathrm{comp}}(n,m)\!\leq\!2\lceil\cK'(n,m)/2\rceil$,
the multiple of $2$ can be removed in the smooth case, and 
$P_{\mathrm{comp}}(1,m)\!=\!m-1$. }
\item{$P_{\mathrm{non}}(n,m)\!\leq\!2(P_{\mathrm{comp}}(n-1,m)
+P_{\mathrm{non}}(n-1,m))$, $P_{\mathrm{non}}(2,m)\!\leq\!\lceil m/2\rceil$, 
and $P_{\mathrm{non}}(1,m)$ is $1$ or $0$ according as $m$ is $0$ or not. } 
\item{$(n,m)\!=\!(2,3) \Longrightarrow Z$ has no more than 
$3$ inflection points and no more than $1$ isolated point of vertical 
tangency.}  
\end{enumerate} 
\end{thm}  
\begin{rem} 
Note that a non-compact component of $Z$ can actually have compact closure, 
since $\R^n_+$ is not closed in $\Rn$, e.g., $\{ (x_1,x_2)\!\in\!\R^2_+ 
\; | \; x^2_1+x^2_2\!=\!1\}$. Also, Bertrand Haas has pointed out  
that the bound on $P_\mathrm{non}(2,m)$, at least in the case of integral 
exponents, may date back to work of Isaac Newton in the 17$^\thth$ century 
on power series. $\diamond$ 
\end{rem}

While the above bounds on the number of connected components are non-explicit, 
they are stated so they can immediately incorporate any advance in 
computing $\cK'(n,m)$. So for a general and explicit 
bound independent of the underlying polynomial degrees now, one could, 
for instance, simply insert the explicit upper bound for $\cK'(n,\mu)$ 
appearing in Khovanski's Theorem on Fewnomials (see section \ref{sub:rel} 
below). 
\begin{cor}
Following the notation of theorem \ref{thm:cool}, 
$Z$ has no more than \scalebox{.85}[1]{$2^{(m-1)(m-2)/2}2^{n-2}n(n+1)^{m-1}$} 
connected components. In particular, a curve in the positive quadrant defined 
by a tetranomial has no more than $4$ compact (resp.\ $2$ non-compact) 
connected components. \qed 
\end{cor} 
\noindent 
The bound above is already significantly sharper than an earlier 
bound of \scalebox{.85}[1]{$2^{m(m-1)/2} (2n)^{n-1}(2n^2-n+1)^m$}, 
which held only for the smooth case, following from \cite[sec.\ 3.14, 
cor.\ 5]{few}.  
The bounds of theorem \ref{thm:cool} 
are further refined in theorem \ref{thm:bounds} of section \ref{sec:morse}, 
and these additional bounds also improve an earlier result of the author 
on smooth algebraic hypersurfaces \cite[cor.\ 3.1]{real}. 

\subsection{Important Related Results}\mbox{}\\ 
\label{sub:rel}  
It is interesting to note that the best current general bounds in the 
direction of GKC are exponential in the number of monomial terms of $F$, 
even for fixed $n$. Observe one of the masterpieces of real algebraic geometry. 
\begin{kho}
(See also \cite{kho} and \cite[cor.\ 7, sec.\ 3.12]{few}.)  
Let $F$ be an $n\times n$ fewnomial system and $\mu$ the 
total number of distinct exponent vectors of $F$. Then $F$ 
has no more than $(n+1)^\mu2^{\mu(\mu-1)/2}$ non-degenerate roots in the 
positive orthant, i.e., $\cK'(n,\mu)\!\leq\!(n+1)^\mu2^{\mu(\mu-1)/2}$. \qed 
\end{kho} 

\begin{rem}
\label{rem:kho} 
In the case $(n,m_1,m_2)\!=\!(2,3,m)$, one can 
divide both equations by suitable monomials to obtain $\mu\!=\!m+2$ 
and thus $\cN'(3,m)\!\leq\!\cK'(3,m+2)$. So Khovanski's bound implies 
$\cN'(3,m)\!\leq\!3^{m+2}2^{(m+2)(m+1)/2}$. It is also very easy to 
see that a simple application of Gaussian elimination yields 
$\cK'(n,m)\!\leq\!\cN'(\underset{n}{\underbrace{m-1,\ldots,m-1}})$.  
$\diamond$ 
\end{rem} 

Non-trivial lower bounds on even $\cN'(3,m)$ are scarce and surprisingly little 
else is known about what an optimal version of Khovanski's Theorem on 
Fewnomials should resemble. For example, an earlier (conjectural) polyhedral 
generalization of Descartes' 
Rule to multivariate systems of equations proposed by Itenberg and Roy 
in 1996 \cite{ir} (based on a famous construction of Oleg Viro from 1989 and 
extensions by Bernd Sturmfels \cite{viro}) was recently disproved 
\cite{liwang}. Also, a bit earlier, Bernd Sturmfels bet (and unfortunately 
lost) US\$500 on a challenge problem involving a family of polynomial systems 
of type $(4,4)$ \cite{lr}. 

To the best of the authors' knowledge, all other general bounds on the number 
of real roots depend strongly on the individual exponents of $F$ and are 
actually geared more toward counting complex roots, e.g., 
\cite{bkk,kaza,blr,jpaa,real}. So 
even proving $\cN(3,3)\!<\!\infty$ already requires a different approach. 
Nevertheless, the aforementioned bounds can be quite practical when the 
exponents are integral and the degrees of the polynomials are small. 

In any event, it still remains unknown whether $\cK'(n,m)$ is 
polynomial in $m$ for $n$ fixed. (The polynomial 
system \mbox{$(x^2_1-3x_1+2,$} $\ldots,x^2_n-3x_n+2)$ shows us that fixing 
$n$ is necessary.) Even the case of a trinomial and an $m$-nomial, in two 
variables, remains open. More to the point, it is also unknown whether a 
simple modification (e.g., increasing the original GKC bound by a 
constant power or a factor exponential in $n$) changes the status of GKC 
from false to true. The $2k\times 2k$ fewnomial system 
$(x^{108}_1+1.1y^{54}_1-1.1y_1,y^{108}_1+1.1x^{54}_1-1.1x_1,\ldots,
x^{108}_k+1.1y^{54}_k-1.1y_k,y^{108}_k+1.1x^{54}_k-1.1x_k)$, thanks to 
Haas' counter-example (cf.\ remark \ref{rem:haas}), easily shows that the 
GKC bound now needs at least an extra multiple no smaller than 
$\left(\frac{\sqrt{5}}{2}\right)^n$ if it is to be salvaged.  
\begin{rem} 
\label{rem:napo} 
Domenico Napoletani has recently shown that to calculate 
$\cN'(m_1,\ldots,m_n)$ for any given $(m_1,\ldots,m_n)$, it suffices 
to restrict to the case of integral exponents \cite{napo}. Here, we 
will bound $\cN(3,m)$ directly without using this reduction. $\diamond$ 
\end{rem} 

Let us also make a related number-theoretic 
observation: Hendrik W.\ Lenstra has shown that for any fixed number field 
$\cL$, the maximal (finite) number of roots in $\cL$ of a univariate 
$m$-nomial, with integral exponents and coefficients in $\cL$, is 
quasi-quadratic in $m$ and 
independent of the degree of the polynomial \cite{lenstra2}. Thus 
an immediately corollary of theorem \ref{thm:tri3} (and theorem \ref{thm:tri1} 
and remark \ref{rem:extra} of section \ref{sec:back}) is that Lenstra's 
result can be effectively extended to certain families of fewnomial systems, 
provided we fix $n$ and restrict to real algebraic number fields. (Fixing $n$ 
is necessary for the same reason as in the last paragraph.) 

Whether Lenstra's result can be more fully 
extended to polynomial systems is also an open question, even in 
the case of two bivariate trinomials. 
However, it is at least now known that 
that the number of {\bf geometrically}\footnote{A root is geometrically 
isolated iff it is a zero-dimensional component of the underlying 
zero set in $\bar{\cL}^n$, where $\bar{\cL}$ is the algebraic closure of 
$\cL$.} isolated roots in 
$\cL^n$ of any $k\times n$ polynomial system can be bounded above by {\bf 
some} function depending only on $\cL$, $n$, and the total number of distinct 
exponent vectors \cite{myadic}.\footnote{ In fact, as was done more explicitly 
in \cite{lenstra2} 
for the univariate case, one can also allow $\cL$ to be any finite extension 
of the $p$-adic rationals. The latter setting is perhaps closer to our 
current focus since $\R$, like the $p$-adics, is a metrically complete field.}  

\subsection{Organization of the Proofs}\mbox{}\\ 
Section \ref{sec:back} provides some background and unites some simple cases 
where GKC in fact holds. 
We then prove theorems \ref{thm:tri3} and \ref{thm:cool} in sections 
\ref{sec:tx} and \ref{sec:morse}, respectively. Proving the upper 
bound on $\cN(3,m)$ turns out to be surprisingly elementary, but 
lowering the bound on $\cN(3,3)$ to $5$ then becomes a more 
involved case by case analysis. 

Section \ref{sec:tri} then gives an alternative geometric proof that 
$\cN(3,3)\!\leq\!6$. We include this second proof for motivational purposes, 
since it 
was essentially the first improvement we found over $\cN'(3,3)\!\leq\!248832$. 
We then derive bounds for the number of isolated singularities and inflection 
points of an  $m$-nomial curve, and discuss how the underlying {\bf Newton 
polygons} (cf.\ the next section) strongly control how $\cN(3,3)$ can exceed 
$4$ (cf.\ corollary \ref{cor:poly} of section \ref{sec:tri}). Roughly 
speaking, we show that if a fewnomial system of type $(3,3)$ has maximally 
many roots in the positive quadrant, then its underlying exponent vectors must 
be in ``general position.'' In particular, just like Haas' counter-example, the 
underlying Newton polygons of any counter-example to this case of GKC must have 
{\bf Minkowski sum} a hexagon (cf.\ sections \ref{sec:back} and 
\ref{sec:tri}). 

\section{The Pyramidal, Simplicial, and Zero Mixed Volume Cases}
Consider the following constructions. 
\label{sec:back} 
\begin{dfn}
For any $S\!\subseteq\!\Rn$, let $\pmb{\conv(S)}$ denote the smallest 
convex set containing $S$. Also, for any $m$-nomial of the form 
$f\!:=\!\sum_{a\in A} c_ax^a$, we call $\pmb{\supp(f)}\!:=\!\{a \; | \; 
c_a\!\neq\!0\}$ the {\bf support} of $f$, and define 
$\pmb{\newt(f)}\!:=\!\conv(\supp(f))$ to be the {\bf Newton polytope} of $f$. 
More generally, a {\bf polytope} is simply the convex hull of any finite  
point set in $\Rn$. \qed 
\end{dfn}  
\begin{dfn}
Let $F\!=\!(f_1,\ldots,f_n)$ be a fewnomial system and for all
$i$ let $L_i$ be the linear subspace affinely generated by $\supp(f_i)$.
We call $F$ {\bf pyramidal} iff the following condition holds  
for all $i$: either $L_i\!\supsetneq\!L_j$ for all $j\!\neq\!i$, or  
there is a $j$ such that $L_i\oplus L_j\!=\!L_i\oplus L$ for 
some line $L\!\not\subseteq\!L_i$ with $\bO\!\in\!L$. Finally, we call any 
change of variables of the form 
$(x_1,\ldots,x_n)\mapsto y^A\!:=\!(y^{a_{11}}_1\cdots 
y^{a_{n1}}_n,\ldots, y^{a_{1n}}_1\cdots y^{a_{nn}}_n)$, with $A\!:=\![a_{ij}]$ 
a real $n\times n$ matrix, a {\bf monomial change of variables}.  \qed
\end{dfn}

For example, the polynomial systems from remark 
\ref{rem:pyr} are all pyramidal, but the systems from remarks  
\ref{rem:degen} and \ref{rem:haas} are not pyramidal (cf.\ section 1). 
Pyramidal systems are a simple generalization of the so-called
``triangular'' systems popular in Gr\"obner-basis papers on
computer algebra. The latter family of systems simply consists of those $F$
for which the equations and variables can be reordered so that 
for all $i$, $f_i$ depends only on $x_1,\ldots,x_i$. Put another 
way, pyramidal systems are simply the image of a triangular system (with {\bf 
real} exponents allowed) after multiplying the individual equations by 
arbitrary monomials, shuffling the equations, and then performing a monomial 
change of variables. In particular, we note the following elementary 
fact on monomial changes of variables.
\begin{prop}
\label{prop:easy} 
If $A$ is a real non-singular $n\times n$ matrix, then 
$(x^{A})^{A^{-1}}\!=\!x$ and the map $x\mapsto x^A$ is an 
analytic automorphism of the positive orthant. In particular, 
such a map preserves smooth points, singular points, and the 
number of compact and non-compact connected components, of 
analytic subvarieties of the positive orthant. Furthmore,  
this invariance also holds for fewnomial zero sets in the 
positive orthant. \qed 
\end{prop} 
\noindent 
The assertion on analytic subvarieties follows easily from an application of 
the chain rule from calculus, and noting that such monomial maps 
are also diffeomorphisms. That the same invariance holds for fewnomial 
zero sets follows immediately upon observing that the substition 
$(x_1,\ldots,x_n\!=\!(e^{z_1},\ldots,e^{z_n})$ maps any $n$-variate real 
$m$-nomial to a real analytic function, and noting that 
$(t_1,\ldots,t_n)\mapsto (e^{t_1},\ldots,e^{t_n})$ is a diffeomorphism from 
$\Rn$ to $\Rn_+$. 
\begin{rem}
The zero set of $x_1+x_2-1$, and the change of variables 
$(x_1,x_2)\mapsto (\frac{y_1}{y_2},y_1y_2)$, show that 
the number of isolated {\bf inflection} 
points need not be preserved by such a map: the underlying curve 
goes from having no isolated inflection points  
to having one in the positive quadrant. $\diamond$ 
\end{rem} 

We will also need the following analogous geometric extension of 
the concept of an over-determined system.  
\begin{dfn} 
\label{dfn:zero} 
Given polytopes $P_1,\ldots,P_n\!\subset\!\Rn$, we say that they have {\bf 
mixed volume\footnote{The reader curious about mixed volumes of polytopes in 
this context can consult \cite{buza,jpaa} for further discussion.} zero} 
iff for some $d\!\in\!\{0,\ldots,n-1\}$ there exists a $d$-dimensional 
subspace of $\Rn$ containing translates of $P_i$ for at least $d+1$ distinct 
$i$. \qed  
\end{dfn} 
\noindent 
A simple special case of an $n$-tuple of polytopes with mixed volume zero 
is the $n$-tuple of Newton polytopes of an $n\times n$ fewnomial system 
where, say, the variable $x_i$ does not appear. Indeed, by multiplying 
the individual $m$-nomials by suitable monomials, and applying a 
suitable monomial change of variables, the following 
corollary of proposition \ref{prop:easy} is immediate. 
\begin{cor} 
\label{cor:zero} 
Suppose $F$ is a fewnomial system, with only finitely many roots 
in the positive orthant, whose $n$-tuple of Newton polytopes has mixed volume 
zero. Then $F$ has no roots in the positive orthant. \qed 
\end{cor} 
\noindent 
Indeed, modulo a suitable monomial change of variables, one need only observe 
that the existence of a single root in the positive orthant implies the 
existence of an entire ray of roots (parallel to some coordinate axis) in the 
positive orthant. 

We will also need the following elegant extension of Descartes' Rule to real 
exponents. It's proof involves a very simple induction using Rolle's 
Theorem (cf.\ the next section) and dividing by suitable monomials 
\cite{few} --- tricks we will build upon in the next section. 
\begin{dfn}
For any sequence $(c_1,\ldots,c_m)\!\in\!\R^m$, it's {\bf number of 
sign alternations} is the number of pairs $\{j,j'\}\!\in\{1,\ldots,m\}$ such 
that $j\!<\!j'$, $c_jc_{j'}\!<\!0$, and $c_i\!=\!0$ when $j\!<\!i\!<\!j'$. \qed 
\end{dfn} 
\begin{des}
\label{thm:des} 
Let $c_1,a_1,\ldots,c_m,a_m$ be any real numbers with $a_1\!<\cdots <\!a_m$. 
Then the number of positive roots of $\sum^m_{i=1}c_ix^{a_i}_1$ is 
at most the number of sign alternations in the sequence 
$(c_1,\ldots,c_m)$. In particular, $\cN'(m)\!=\!\cN(m)\!=\!m-1$. \qed 
\end{des} 

As a warm-up, we can now prove a stronger version of GKC for the following 
families of special cases.  
\begin{thm}
\label{thm:tri1} 
Suppose $F$ is an $n\times n$ fewnomial system of type 
$(m_1,\ldots,m_n)$ (so $m_1,\ldots,m_n\!\geq\!1$) and we restrict to those 
$F$ which also satisfy one of the following conditions: 
\begin{itemize} 
\item[{\bf (a)}]{The $n$-tuple of Newton polytopes of $F$ 
has mixed volume zero.} 
\item[{\bf (b)}]{All the supports of $F$ can be translated into 
a single set of cardinality $\leq\!n+1$.}  
\item[{\bf (c)}]{$F$ is pyramidal.} 
\end{itemize} 
Then, following the notation of theorem \ref{thm:tri3}... 
\begin{enumerate} 
\addtocounter{enumi}{-1} 
\item{$\cN(m_1,\ldots,m_n)$ is respectively $0$, $1$, or 
$\prod^n_{i=1}(m_i-1)$ in case (a), (b), or (c).} 
\item{In cases (a), (b), and (c), $F$ has infinitely many roots 
$\Longrightarrow F$ has {\bf no} isolated roots.}
\item{In general, $\cN'(1,m_2,\ldots,m_n)\!=\!\cN(1,m_2,\ldots,m_n)\!=\!0$, 
$\cN'(2,m_2,\ldots,m_n)\!=\!\cN'(m_2,\ldots,m_n)$, 
and $\cN(2,m_2,\ldots,m_n)\!=\!\cN(m_2,\ldots,m_n)$.}  
\item{$\cN'(2,m)\!=\!\cN(2,m)\!=\!m-1$.}  
\end{enumerate} 
\end{thm} 

\noindent 
{\bf Proof:} First note that the Newton polytopes must all be nonempty. 
The case (a) portion of assertions (0) and (1) then follows immediately 
from corollary \ref{cor:zero}. Note also that the case (a) portion 
of assertion (0) immediately implies our formula for $\cN(1,m_2,\ldots,m_n)$ 
(and thus $\cN'(1,m_2,\ldots,m_n)$ as well) in assertion (2), since the 
underlying $n$-tuple of polytopes clearly has mixed volume zero. 

The case (b) portion of assertions (0) and (1) follows easily upon observing 
that $F$ is a linear system of $n$ equations in $n$ monomial terms, after 
multiplying the individual equations by suitable monomial terms. We can then 
finish by proposition \ref{prop:easy}. 

To prove the case (c) portion of assertions (0) and (1), note that the case 
$n\!=\!1$ follows immediately from UGDRS. For $n\!>\!1$, we have the following 
simple proof by induction: Assuming GKC holds for all $(n-1)\times (n-1)$ 
pyramidal systems, consider any $n\times n$ pyramidal system $F$. Then, via a 
suitable   monomial change of variables, multiplying the individual equations 
by suitable monomials, and possibly reordering the $f_i$, we can assume that 
$f_1$ depends only on $x_1$. (Otherwise, $F$ wouldn't be pyramidal.) We thus 
obtain by UGDRS that $f_1$ has at most $m_1-1$ positive roots. By 
back-substituting these roots into $F'\!:=\!(f_2,\ldots,f_n)$, we obtain a 
new $(n'-1)\times (n'-1)$ pyramidal fewnomial system of type 
$(m'_2,\ldots,m'_{n'})$ with $n'\!\leq\!n$ and 
$m'_2\!\leq\!m_2,\ldots,m'_n\!\leq\!m_n$. 
By our induction hypothesis, we obtain that each such specialized $F'$ has at 
most $\prod^{n'}_{i=2}(m'_i-1)$ isolated roots in the positive orthant, and 
thus $F$ has at most $\prod^n_{i=1}(m_i-1)$ isolated roots in the positive 
orthant. (Remark \ref{rem:pyr} from the introduction shows us that this bound 
can indeed be attained.) 

Our recursive formulae for $\cN'(2,m_2,\ldots,m_n)$ and 
$\cN(2,m_2,\ldots,m_n)$ from assertion (2) then follow by applying just the 
first step of the preceding induction argument, and noting that proposition 
\ref{prop:easy} tells us that our change of variables preserves non-degenerate 
roots. 

Assertion (3) follows immediately from assertion (2) via UGDRS. \qed 

\begin{rem} 
\label{rem:extra} 
One can of course combine and interweave families (a), (b), and (c) 
to obtain less trivial examples where GKC is true. More generally, 
one can combine theorems \ref{thm:tri3} and \ref{thm:tri1} 
to obtain bounds significantly sharper than Khovanski's Theorem on 
Real Fewnomials, free from Jacobian assumptions, for additional families of 
fewnomial systems. $\diamond$ 
\end{rem} 

\section{Substitutions and Calculus: Proving 
Theorem \ref{thm:tri3}} 
\label{sec:tx} 
Let us preface our first main proof with some useful basic results. 
\begin{lemma}
\label{lemma:tri2} 
For $m_1\!=\!1+\dim \newt(f_1)$, the computation of $\cN'(m_1,\ldots,m_n)$ 
and $\cN(m_1,\ldots,m_n)$ can be reduced to the case where 
$f_1:= 1\pm x_1\pm \cdots \pm x_{m_1-1}$ (with the signs in $f_1$ {\bf not} all 
``$+$'') and, for all $i$, $f_i$ has $1$ as one of its monomial terms. 
In particular, for $m_1\!=\!3$, we can assume further that 
$f_1\!:=\!1-x_1-x_2$. 
\end{lemma} 
 
\noindent 
{\bf Proof:} By dividing each $m_i$-nomial by a suitable 
monomial term, we can immediately assume that 
all the $f_i$ possess the monomial term $1$. 
In particular, we can also assume that the origin $\bO$ is a vertex of 
$\newt(f_1)$. Note also that the sign condition on $f_1$ must obviously hold, 
for otherwise the value of $f_1$ would be positive on the positive 
orthant. (The refinement for $m\!=\!3$ then follows by picking the monomial 
term one divides $f_1$ by a bit more carefully.) So we now need only check 
that the desired canonical form for $f_1$ can be attained. 

Suppose $f_1\!:=\!1+c_1x^{a_1}+\cdots +c_{m_1-1}x^{a_{m_1-1}}$. 
By assumption, $\newt(f_1)$ is an $m_1$-simplex with vertex set 
$\{\bO,a_1,\ldots,a_{m_1-1}\}$, so $a_1,\ldots,a_{m_1-1}$ 
are linearly independent. Now pick any $a_{m_1},\ldots,a_n\!\in\!\Rn$ 
so that $a_1,\ldots,a_n$ are linearly 
independent. The substitution $x\mapsto x^{A^{-1}}$ (with $A$ the matrix whose 
columns are 
$a_1,\ldots,a_n$) then clearly sends $f_1 \mapsto 1+c_1 x_1+\cdots 
+c_{m_1-1} x_{m_1-1}$, and proposition \ref{prop:easy} tells us that 
this change of variables preserves degenerate and non-degenerate roots  
in the positive orthant. Then, via the change of variables 
$(x_1,\ldots,x_{m_1-1}) \mapsto (x_1/|c_1|,\ldots,x_{m_1-1}/|c_{m_1-1}|)$, we 
obtain that $f_1$ can indeed be placed in the desired form. (The latter change 
of variables preserves degenerate and non-degenerate roots in the positive 
orthant for even more obvious reasons.) \qed 

Recall that a polynomial $p\!\in\!\C[x_1,\ldots,x_n]$ is {\bf homogeneous} 
iff $p(ax_1,\ldots,ax_n)\!=\!a^dp(x_1,\ldots,x_n)$ for some non-negative 
integer $d$. 
\begin{prop} 
\label{prop:ind} 
Suppose $p\!\in\!\R[x_1,x_2]$ is homogeneous and has degree $d\!\geq\!0$. 
Also let $\alpha,\beta\!\in\!\R$. 
Then there is a homogeneous $q\!\in\!\R[x_1,x_2]$, either identically zero or 
of degree $d+1$, such that $\frac{d}{dt}\left(t^\alpha (1-t)^\beta 
p(t,1-t)\right)\!=\! t^{\alpha-1}(1-t)^{\beta-1}q(t,1-t)$. In particular, 
$\frac{d}{dt}\left(t^\alpha (1-t)^\beta p(t,1-t)\right)$ is identically 
zero iff $p(x_1,x_2)\!=\!x^{-\beta}_1x^\alpha_2$, $\alpha-\beta\!=\!d$, 
and $\alpha,\beta\!\in\!\Z$ with $\alpha\!\geq\!0$. \qed 
\end{prop} 

\noindent 
{\bf Proof:} By the chain-rule, $\frac{d}{dt}\left(t^\alpha 
(1-t)^\beta p(t,1-t)\right)$ is simply \[\alpha t^{\alpha-1}(1-t)^\beta p(t,1-t) 
+ \beta t^\alpha (1-t)^{\beta-1}p(t,1-t) + t^\alpha (1-t)^\beta 
(p_1(t,1-t)-p_2(t,1-t)),\] 
where $p_i$ denotes the partial derivative 
of $p$ with respect to $x_i$. Factoring out a multiple of 
\mbox{$t^{\alpha-1}(1-t)^{\beta-1}$} from the preceding expression, we then 
easily obtain that we can in fact take 
\[q(x_1,x_2)\!=\!(\alpha x_2+\beta x_1)p(x_1,x_2)+x_1x_2(p_1(x_1,x_2)-p_2(x_1,x_2)). \] 
The final assertion of our proposition then follows immediately. \qed  

\begin{rolle}\footnote{For a simple proof, note that the special case 
$r\!=\!2$ follows immediately from the Mean Value Theorem of calculus 
(see, e.g., \cite[thm.\ 5.10, pg.\ 107]{rudin}), 
since we can replace $[a,b]$ by a sub-interval whose end-points are 
the roots of $g$. The general case then follows by replacing $a$ (resp.\ 
$b$) by the smallest (resp.\ largest) root of $g$, and then subdividing 
$[a,b]$ into $r-1$ sub-intervals whose endpoints consist of the roots of $g$.}  
Let $g : [a,b] \longrightarrow \R$ be any function with a well-defined  
derivative $g'$ defined on $(a,b)$. Then $g$ has $r$ roots in 
$[a,b] \Longrightarrow g'$ has at least $r-1$ roots in $(a,b)$. \qed 
\end{rolle}
\begin{lemma}
Let $k\!\geq\!2$. 
Then for any real $c_1,a_1,b_1,\ldots,c_k,a_k,b_k$, the function
\[f(t):=1+c_1 t^{a_1} (1-t)^{b_1}+\cdots + c_k t^{a_k}(1-t)^{b_k}\]
has at most $2^{k+1}-2$ roots in the open interval $(0,1)$. Furthermore, 
$f$ has exactly $r$ roots in $(0,1) \Longrightarrow$ there exist 
$\tilde{c}_1,\ldots,\tilde{c}_k\!\in\!\R$ such that 
\[\tilde{f}(t)\!:=\!1+\tilde{c}_1 t^{a_1} (1-t)^{b_1}+\cdots + 
\tilde{c}_k t^{a_k}(1-t)^{b_k}\] 
has at least $r$ roots in $(0,1)$, and no root of $\tilde{f}$ is degenerate.  
\label{lemma:ind}
\end{lemma}

\noindent 
{\bf Proof:} 
Henceforth, let us assume all roots lie in the open interval $(0,1)$. 
Assume $f$ has exactly $r$ roots. Then by Rolle's Theorem, $f'$ has at 
least $r-1$ roots. Since \[f'(t)= \sum^k_{i=1}c_i t^{a_i-1} (1-t)^{b_i-1}
(a_i(1-t)+b_it),\] and since $t^{a_k-1} (1-t)^{b_k-1}$ never vanishes in 
$(0,1)$, the function 
\[ g_1(t)\!:=\!c_k(a_k(1-t)+b_kt)+\sum^{k-1}_{i=1}c_i t^{a_i-a_k} 
(1-t)^{b_i-b_k}(a_i(1-t)+b_it) \] 
has at least $r-1$ roots. 

By Rolle's Theorem again, $g''_1$ has at least $r-3$ roots.
By proposition \ref{prop:ind}, $g''_1(t)$ will then be of the form
$\sum^{k-1}_{i=1}c_i t^{a_i-a_k-2} (1-t)^{b_i-b_i-2}q_{i,1}(t,1-t)$, 
where the $q_{i,1}\!\in\!\R[x_1,x_2]$ are homogeneous polynomials, which are 
either identically zero or of degree $3$. In 
particular, we can assume that at least one $q_{i,1}$ must be different from 
the zero polynomial. (For otherwise we would obtain that 
$g''_1\!=\!0$ identically, which would in turn imply that $g_1$ is a linear 
function, and thus $r\!\leq\!3\!<\!2^{k+2}-2$.)  
By again dividing by a suitable monomial in $t$ and $1-t$, we then see that 
$g''_1$ has the same number of roots as
\[ g_2(t):= q_{k-1,1}(t,1-t)+\sum^{k-2}_{i=1}c_i t^{a_i-a_{k-1}} 
(1-t)^{b_i-b_{k-1}}q_{i,1}(t,1-t).\] 
Thus $g_2$ has at least $r-3$ roots. 

By induction, we then easily obtain a sequence of 
polynomials $g_1,\ldots,g_{j}$, where $j\!\leq\!k$ and 
$g_j\!=\!q_{1,j-1}(t,1-t)$ for some homogeneous $q\!\in\!\R[x_1,x_2]$ 
of degree $2^j-1$ having at least $r-(2^j-1)$ roots.  So by Rolle's 
Theorem one last time, $r\!\leq\!(2^j-1)+(2^j-1)\!\leq\!2^{k+1}-2$ and we 
are done with the first part of our lemma. 

To prove the second part, note that the first part of our lemma implies 
that $f$ has only finitely many critical values --- no more than $2^{k+1}-2$, 
in fact. So for all $\delta\!\in\!\Rs$ with $|\delta|$ 
sufficiently small, $f-\delta$ will have {\bf no} degenerate roots. We can 
in fact guarantee that $f-\delta$ will also have at least $r$ non-degenerate 
roots in $(0,1)$ as follows: Let $n_+$ 
(resp.\ $n_-$) be the number of roots $t$ of $f$ with $f'(t)\!=\!0$ 
and $f''(t)\!>\!0$ (resp.\ $f''(t)\!<\!0$). Clearly then, 
for all $\delta\!\in\!\Rs$ with $|\delta|$ sufficiently small, 
$f-\delta$ will have exactly $r+n_--n_+$ or $r+n_+-n_-$ roots, according as 
$\delta\!>\!0$ or $\delta\!<\!0$. (The analogous statement for roots in 
$(0,1)$ holds as well, since $(0,1)$ is open.) So let $\tilde{\delta}$ be 
sufficiently small, and of the correct sign, so that $f-\tilde{\delta}$ has at 
least $r$ roots in $(0,1)$ and {\bf no} degenerate roots.

To conclude, simply let $\tilde{c}_i\!:=\!\frac{c_i}{1-\tilde{\delta}}$ for 
all $i$.  Since $\tilde{f}$ is thus $\frac{f-\tilde{\delta}}{1-
\tilde{\delta}}$, we are done. \qed 

\medskip 

\noindent 
{\bf Proof of Theorem \ref{thm:tri3}:} 
First note that by lemma \ref{lemma:tri2}, we can immediately 
reduce to the case of a fewnomial system of the form $F\!:=\!(1 - x_1 - x_2, 
1+c_1 x_1^{a_1} x_2^{b_1} +\cdots+c_{m-1}x_1^{a_{m-1}}x^{b_{m-1}}_2)$, 
and this reduction preserves the degeneracy or non-degeneracy of any 
root of $F$. We can then simply solve for $x_2$ via the first 
equation and then substitute into the second equation to obtain a bijection 
between the roots of $F$ in the positive quadrant and the roots of 
$f(t):=1+c_1 t^{a_1}(1-t)^{b_1}+\cdots+c_{m-1}t^{a_{m-1}} (1-t)^{b_{m-1}}$ with 
$0<t<1$. A simple Jacobian calculation yields that 
$(\zeta_1,\zeta_2)$ is a degenerate root of $F \Longleftrightarrow$  
[$\sum^{m-1}_{i=1}c_i\zeta^{a_i-1}_1(1-\zeta_1)^{b_i-1} 
(a_i(1-\zeta_1)-b_i\zeta_1)\!=\!0$ and $f(\zeta_1)\!=\!0$] 
$\Longleftrightarrow f'(\zeta_1)\!=\!f(\zeta_1)\!=\!0$. So degenerate roots 
of our univariate reduction correspond bijectively to degenerate roots of $F$.  

By lemma \ref{lemma:ind}, and the fact that $\cN'(3,m)\!\leq\!\cN(3,m)$, we 
immediately obtain $\cN(3,m)\!\leq\!2^m-2$ and $\cN'(3,m)\!=\!\cN(3,m)$.  
Our upper bounds on $\cN(3,4)$ and $\cN(3,5)$ are then simply specializations 
of our new upper bound for $\cN(3,m)$. 

To now prove that $\cN(3,3)\!=\!5$, thanks to Haas' counter-example from 
remark \ref{rem:haas}, it suffices to show that $\cN(3,3)\!<\!6$. 
To do this, let us specialize our preceding notation to 
$m\!=\!3$, $(c_1,c_2)\!=\!(-A,-B)$, and $(a_1,b_1,a_2,b_2)\!=\!(a,b,c,d)$, 
for some $a,b,c,d\!\in\!\R$ and positive $A$ and $B$. (Restricting to 
positive $A$ and $B$ can easily be done simply by dividing $f_2$ by a 
suitable monomial term, \`a l\'a the proof of lemma \ref{lemma:tri2}.) 

By using symmetry we can then clearly reduce to the following cases:\\ 
\begin{tabular}{cc}
{\bf A.} $a,b,c>0$ and $d<0$ & {\bf B.} $a,c>0$ and 
$b,d<0$\hspace{1.85cm}\mbox{}\\
{\bf C.} $a,b>0$ and $c,d<0$ & {\bf D.} $a,b,c,d>0$\hspace{3.1cm}\mbox{}\\
{\bf E.} $a>0$ and $b,c,d<0$ & {\bf F.} $a,b,c,d<0$\hspace{3.1cm}\mbox{}\\
{\bf G.} $a,d>0$, $b,c<0$\hspace{.6cm}\mbox{} & \mbox{}\hspace{1.9cm}{\bf H.} 
At least one of the numbers $a,b,c,d$ is zero. 
\end{tabular}

In particular, our earlier substitution trick tells us that 
it suffices to show that any \[f(t):=1-A t^a(1-t)^b-Bt^c (1-t)^d, \] 
{\bf with all roots non-degenerate}, always has strictly less than $6$ roots 
in the open interval $(0,1)$. So let $r$ be the number of roots of any 
such non-degenerate $f$ in $(0,1)$. 

Let us now prove $r\!<\!6$ in all 8 cases: 
\begin{itemize}
\item[{\bf A.}] $a,b,c>0$, $d<0$:  

Let $Q(x)=1-Ax^a(1-x)^b$ and $R(x)=Bx^c(1-x)^d$. 
The roots of $f$ may be regarded as the intersections of 
$y\!=\!Q(x)$ and $y\!=\!R(x)$ in the positive quadrant.  Since 
$\lim_{x\rightarrow 0^+} Q(x)=1$, 
$\lim_{x\rightarrow 1^-} Q(x)= 1$,   $\lim_{x\rightarrow 0^+} R(x)=0$, and 
$\lim_{x\rightarrow 1^-} R(x)=\infty$, it is easy to see via 
the Intermediate Value Theorem of calculus that the number of 
intersections must be odd. (One need only note that $f\!=\!Q-R$ and 
that the signs of $f'$ at the ordered roots of $f$ alternate.) So $r\!<\!6$.  
\item[{\bf B.}] $a,c>0$, $b,d<0$: \\
By an argument similar to that of case A, $r$ is odd and thus less than $6$.
\item[{\bf C.}] $a,b>0$, $c,d<0$: \\See lemma \ref{l1} below.
\item[{\bf D.}] $a,b,c,d>0$: \\See lemma \ref{l2} below.
\item[{\bf E.}] $a>0$, $b,c,d<0$:\\See lemma \ref{l3} below.
\item[{\bf F.}] $a,b,c,d<0$. \\ Multiplying $f(t)$ by 
$t^{\max\{-a,-c\}}(1-t)^{\max\{-b,-d\}}$, we can immediately reduce to 
case D.  
\item[{\bf G.}] $a,d>0$, $b,c<0$:\\See lemma \ref{l4} below.
\item[{\bf H.}] At least one of the numbers $a,b,c,d$ is zero:\\
Use lemma \ref{lemma:basic} below, noting that our hypotheses 
here imply that either $F$ or $\hat{F}$ is a 
quadratic polynomial. 
\end{itemize}
This concludes the proof of theorem \ref{thm:tri3}. \qed

We now detail the lemmata we cited above. 

\begin{lemma}
\label{lemma:basic} 
Following the notation of the proof of theorem \ref{thm:tri3}, recall 
that $\pmb{r}$ is the number of roots of\\ $f(t)\!:=\!1-At^a(1-t)^b-Bt^c(1-t)^d$ in 
the open interval $(0,1)$, where $f$ has no degenerate roots. Also let 
$g(t):=\frac{A}{B} t^{a-c}(1-t)^{b-d}(-a(1-t)+bt)-c(1-t)+dt$,\\
\scalebox{.75}[1]{$F(u):=  -a(a-c)(a-c-1)u^3 
+(a-c) [2a(b-d+1)+b(a-c+1)]u^2+(d-b)[a(b-d+1)+2b(a-c+1)]u +b(b-d)(b-d-1)$}, 
and\\ 
\scalebox{.75}[1]
{$\hat{F}(u):= -c(c-a)(c-a-1)u^3 +(c-a) [2c(d-b+1)+d(c-a+1)]u^2 
+(b-d)[c(d-b+1)+2d(c-a+1)]u +d(d-b)(d-b-1)$}.  
Finally, let $N$ (resp.\ $M$) be the number of roots in $(0,1)$ of $g$  
(resp.\ the maximum of the number of positive roots of $F$ and  
$\hat{F}$). Then $r-3\!\le\!N-2\!\le\!M\!\leq\!3$. 
\end{lemma} 

\noindent 
{\bf Proof:} Just as in the proof of lemma \ref{lemma:ind},  
we easily see by Rolle's Theorem and division by suitable 
monomials in $t$ and $1-t$ that $r-1$ is no more than  
the number of roots in $(0,1)$ of $g$. So $r-1\!\leq\!N$. 
Note also that, in a similar way, $r-1$ is no more than the number 
of roots of $\hat{g}(t)\!:=\!\frac{B}{A}t^{c-a}(1-t)^{d-b}g(t)$ in 
$(0,1)$, and the latter function has the same number of roots in 
$(0,1)$ as $g$. 

To conclude, simply note that for suitable 
$\alpha,\beta,\gamma,\delta\!\in\!\R$, we have that 
$F(\frac{1-t}{t})\!=\!t^\alpha(1-t)^\beta g''(t)$ and 
$\hat{F}(\frac{1-t}{t})\!=\!t^\gamma(1-t)^\delta \hat{g}''(t)$. 
So, by our preceding trick again, $N-2\!\leq\!M$, and thus 
$r-3\!\leq\!M$. That $M\!\leq\!3$ is clear from the fundamental 
theorem of algebra. \qed  

\begin{lemma}
Following the notation of lemma \ref{lemma:basic}, let  
\[T(x):=\frac{A}{B} x^{a-c}(1-x)^{b-d}(-a(1-x)+bx), \qquad S(x):=c-(c+d)x,\]
and
\[\hat{T}(x):=\frac{B}{A} x^{c-a}(1-x)^{d-b}(-c(1-x)+dx),\qquad 
\hat{S}(x)=a-(a+b)x.\] 
Then [$a,b>0$ and $c,d<0$] $\Longrightarrow r<6$. 
\label{l1}
\end{lemma}

\noindent 
{\bf Proof:} By lemma \ref{lemma:basic}, we are done if $M\!<\!3$ or 
$N\!<\!5$. So let us assume $M\!=\!3$ to derive a contradiction. By Descartes' 
Rule of Signs (see section \ref{sec:back} for a generalization), the 
coefficients of $F(u)$ or $\hat{F}(u)$ (ordered by exponent) must have 
alternating signs.  Thus, since $a,a-c,b,b-d\!>\!0$, we have that $a-c-1$ and 
$b-d-1$ must have the same sign. We then need to discuss two cases: 

\begin{itemize}
\item $a-c-1<0$ and $b-d-1<0$: 

This implies $c-a+1>0$ and $d-b+1>0$. Consequently,
coefficients of $u^3$ and $u^2$ in $\hat{F}(u)$ and $F(u)$ are all positive 
--- a contradiction.  
\item $a-c-1>0$ and $b-d-1>0$: 

The roots of $g$ in $(0,1)$ can be regarded as intersections of 
$y\!=\!T(x)$ and  $y\!=\!S(x)$, for $ 0<x<1$.
Since $T(x)\!<\!0$ for $0\!<x\!\ll\!1$ and $T(x)\!>\!0$ for $0\!<1-x\!\ll\!1$, 
there is a smallest positive local minimum $c_0$ of $T$ with 
$T(c_0)\!<0\!$. Thus for $x$ near $c_0$, $T''(x)>0$. Since $T''(x)<0$ for 
$0<x\ll 1$, there is $c^*\!\in\!(0,c_0)$ such that $T''(c^*)=0$.
Let $(x_1,y_1),\ldots,(x_K,y_K)$ be the intersection points
of $y=T(x)$ and $y=S(x)$
with $x_1<x_2<\cdots <x_K$, where a tangent point is counted twice. 
Then for all $i\!\in\!\{1,\ldots,K-1\}$ there is a $c_i\!\in\!(x_i,x_{i+1})$ 
with $T'(c_i)=-(c+d)>0$, and for all $i\!\in\!\{1,\ldots,K-2\}$ there 
is a $d_i\!\in\!(c_i,c_{i+1})$ with $T''(d_i)\!=\!0$. 
Note that $c_0<c_1$. Thus $c^*<d_1$ and therefore 
 $T''(x)=0$ has at least $K-1$ solutions. 
 Since $T''$ and $F$ have the same number of 
positive roots (observing that $T''(u)/F(u)$ is a monomial in $u$ and 
$1-u$), we have $N-1 \le K-1  \le 3 $. \qed 
\end{itemize}
\begin{lemma} Following the notation of lemma \ref{l1},  
$a,b,c,d>0 \Longrightarrow r<6$. 
\label{l2}
\end{lemma}
{\bf Proof:} Again, by lemma \ref{lemma:basic}, we need only show 
that $M\!<\!3$ or $N\!<\!5$. So let us assume $M\!=\!3$. Then by Descartes' 
Rule of Signs, $(a-c)(a-c-1)$ and $(b-d)(b-d-1)$ in the coefficients of $u^3$ 
and $u^0$ in $F(u)$ must have the same sign. There are now four cases to be 
examined.
\begin{itemize}
\item The signs of $a-c$, $a-c-1$, $b-d$, and $b-d-1$ are respectively 
$+,-,+$, and $-$: 

This makes the signs of coefficients of $u^3$ and $u^2$ of
$F(u)$ both positive. 

\item  The signs of $a-c$, $a-c-1$, $b-d$, and $b-d-1$ are respectively 
$-,-,+$, and $+$: 

Since  $b-d>0$, we have $d-b<0$ and $d-b-1<0$. This makes the constant term of
$\hat{F}(u)$ positive, and hence, the coefficients of $u$ and $u^2$ of 
$\hat{F}(u)$ must respectively be negative and positive. That is,
$c(d-b+1)+2d(c-a+1)<0\quad {\rm and} \quad 2c(d-b+1)+d(c-a+1)>0$. 
Thus, $-c(d-b+1)+d(c-a+1)<0.$ This is false, since $b-d-1>0$ and $a-c-1<0.$

\item The signs of $a-c$, $a-c-1$, $b-d$, and $b-d-1$ are all negative: 

By Descartes' rule of signs, $d-b-1$ and $c-a-1$ 
in the coefficients of $y^3$ and $y^0$ of $\hat{F}(y)$  
 must have the same sign. If both are negative, then coefficients of 
$u^3$ and $u^2$ of $F(u)$ would both be negative.
Thus  $d-b-1>0$ and $c-a-1>0$.
It is easy to see that $\hat{T}(x)\!<\!0$ for $0\!<x\!\ll\!1$ and
$\hat{T}(x)\!>\!0$ for $0\!<1-x\!\ll\!1$ and $\lim_{x\rightarrow 0^+} 
\hat{T}(x)=\lim_{x\rightarrow 1^-}\hat{T}(x)=0$.
Now let  $L_0=\min \{c \; | \; 1>c>0, \ \hat{T}(c)<0 \text{ \ 
and \ } c \text{ \ is \ a \ local \ minimum} \}$ and  
$U_0=\max\{c \; | \; 1>c>L_0, \ \hat{T}(c) \text{ \ 
is \ a  \ local \ maximum} \}$. 
Then for $x$ near $L_0$, $\hat{T}''(x)>0$. 
 Since $\hat{T}''(x)<0$ for $0<x\ll 1$, 
 there exists $0<L_1<L_0$ such that $\hat{T}''(L_1)=0$.
Similarly, there is a $U_1\!\in\!(U_0,1)$ such that $\hat{T}''(U_1)\!=\!0$.

The roots of $\frac{B}{A}t^{c-a}(1-t)^{d-b}g$ can be regarded as the 
intersections of $y=\hat{T}(x)$ and  $y=\hat{S}(x)$, for $ 0<x<1$.
Let $(x_1,y_1),\ldots,(x_k,y_k)$ be the intersection points with 
$x_1<x_2<\cdots <x_k$, where a tangent point is counted twice.
Then there exist $ x_i<c_i<x_{i+1}$ such that $\hat{T}'(c_i)=-(a+b)<0$, 
$i=1,\ldots,k-1$
and $ c_i<d_i<c_{i+1}$ such that $\hat{T}''(d_i)=0$, $i=1,\ldots,k-2$.
If $x_1\!>\!L_0$, then $L_1<d_1$. If $x_1<L_0$, then $T(x_1)<0$. This implies 
$T(x_i)<0$ for all $i=1,\ldots,k-2$, since the slope $-(a+b)$ of $\hat{S}(x)$ 
is negative. Therefore, $x_{k-2}< U_0$ and hence $d_{k-2}<U_1$.  
So $\hat{T}''(x)=0$  has at least $k-1$ solutions. 
Since $\hat{T}''(x)=0$ and $\hat{F}(y)=0$ have the same number of 
solutions, we have $N-1 \le k-1 \le M=3$. 

\item The signs of $a-c$, $a-c-1$, $b-d$, and $b-d-1$ are all positive: 
  
Since $a-c-1>0$ and $b-d-1>0$, the proof follows the same line of
arguments as the last case by considering the intersections of $T(x)$ and 
$S(x)$ instead.
\end{itemize}
\hfill \qed

\begin{lemma} Following the notation of lemma \ref{l1}, [$a>0$ and 
$b,c,d<0$] $\Longrightarrow r<6$.
\label{l3}
\end{lemma}
{\bf Proof:}  
Once again, by lemma \ref{lemma:basic}, it suffices to show that 
$M\!<\!3$ or $N\!<\!5$. So let us assume that $M\!=\!3$. 
By checking coefficients of $u^3$ and $u^0$ in $F(u)$, Descartes' Rule 
of Signs tells us that $a-c-1$ and $(b-d)(b-d-1)$ must have different signs. 
There are now three cases to be examined.
\begin{itemize}
\item  $a-c-1$, $b-d$, and $b-d-1$ are all negative.
  
Then the signs of the coefficients of both $u^3$ and $u^2$ in  
$\hat{F}(u)$ will all be positive.
 
\item  The signs of $a-c-1$, $b-d$, and $b-d-1$ are respectively 
$-,+$, and $+$. 

Multiplying $f$ by $x^{-c}(1-x)^{-d}$ yields
$u(x):= x^{-c}(1-x)^{-d}-Ax^{a-c}(1-x)^{b-d}-B$,  
where $-c>0$, $a-c>0$, $-d>1$, and $-d+b>1$. 
The roots of $u$ in $(0,1)$ can be regarded as the intersections of 
the curves 
$y=v(x)=x^{-c}(1-x)^{-d}-Ax^{a-c}(1-x)^{b-d}$ and $y=B$.
Let $(x_1,y_1),\ldots,(x_n,y_n)$ be the intersection points of
$y=v(x)$ and $y=B$ with $x_1<x_2<\cdots <x_n$, where a tangent point is 
counted twice.
Then there exist $ x_i<c_i<x_{i+1}$ such that $v'(c_i)=0$, 
$i=1,\ldots,n-1=\hat{N}$.
Thus $v'$ has at least $\hat{N}$ roots in $(0,1)$. A straightforward 
computation then yields,
\[v'(x)=Ax^{a-c-1}(1-x)^{b-d-1}(-(a-c)(1-x)+(b-d)x)
+x^{-c-1}(1-x)^{-d-1}(-c(1-x)+dx),\]
which clearly has the same number of roots in $(0,1)$ as
\[t(x):=Ax^{a}(1-x)^{b}(-(a-c)(1-x)+(b-d)x)-c(1-x)+dx.\]
Thus $t''$ has at least $\hat{N}-2$ roots in $(0,1)$. Since
\begin{equation*}
\begin{array}{ll}
t'' (x)/A=&x^{a-2}(1-x)^{b-2}[-(a-c)a(a-1)(1-x)^3 
 +a((a+1)(b-d)+2(b+1)(a-c))x(1-x)^2\\
 & -b((b+1)(a-c)+2(b-d)(a+1))x^2(1-x)+(b-d)b(b-1) x^3], 
\end{array}
\end{equation*}
$t''$ has as many roots in $(0,1)$ as
\begin{equation*}
\begin{array}{ll}
P(u)=&-(a-c)a(a-1)u^3 +a((a+1)(b-d)+2(b+1)(a-c))u^2\\
    & -b((b+1)(a-c)+2(b-d)(a+1))u+(b-d)b(b-1) 
\end{array}
\end{equation*}
has positive roots. 
Since $a-1<a-c-1<0$, the coefficients of $u^3$ and $u^0$
in $P(u)$ are both positive. Thus $P$ has at most $2$ positive roots 
and we obtain $\hat{N}-2\le 2$.

\item The signs of $a-c-1$, $b-d$, and $b-d-1$ are respectively $+,+$, and 
$-$: 
  
Since $a-c-1>0$ and $b-d>0$, it is easy to see that $T(x)\!<\!0$ for 
$0\!<\!x\!\ll\!1$ and
$\lim_{x\rightarrow 1^-}T(x)=-\infty$. If $T(x)$ has no local minimum, 
then $y=T(x)$ and $y=S(x)$ have at most one intersection point. Otherwise, 
let $c_0=\min \{c \; | \; 1\!>\!c\!>\!0, \ c \text{ \ is \ a \ local \ 
minimum \ of \ } T\}$. The rest of the proof is similar to that of lemma 
\ref{l1}. \qed
\end{itemize}
\begin{lemma} Following the notation of lemma \ref{l1}, 
[$a,d>0$ and $b,c<0$] $\Longrightarrow r<6$.
\label{l4}
\end{lemma}
{\bf Proof:}
One last time, lemma \ref{lemma:basic} tells us that it suffices to 
prove that $M\!<\!3$ or $N\!<\!5$. So let's assume that $M\!=\!3$. Checking 
signs of coefficients of $u^3$ and 
$u^0$ of both $F(u)$ and $\hat{F}(u)$, Descartes' Rule of Signs tells us that 
$a-c-1<0$ and $d-b-1<0$. On the other hand, the alternating signs of 
coefficients of $u^2$ and $u^1$ of $F(u)$ yield 
$$2a(b-d+1)+b(a-c+1)<0 \quad {\rm and}\quad a(b-d+1)+2b(a-c+1)>0.$$ Thus,
$$-a(b-d+1)+b(a-c+1)=a(d-1)+b(1-c)>0.$$ This is impossible, since 
$d-1<d-1-b<0$, $1-c>0$, $a>0$,  and $b<0$. \qed

\begin{rem} 
When $A=1.12$, $B=0.71$, $a=0.5$, $b=0.02$, 
$c=-0.05$, and $d=1.8$, 
\begin{equation}
f(x)=1-Ax^a(1-x)^b-Bx^c(1-x)^d=0, 0<x<1 
\end{equation}
has $5$ solutions. They are, \mbox{approximately, 
$\{0.00396494, 0.02986317,0.4354707,0.72522344,0.99620026\}$. $\diamond$} 
\end{rem}

\section{A Simple Geometric Approach}
\label{sec:tri}

Let us begin with an extension of Rolle's Theorem to smooth curves
in the plane. 
\begin{lemma} 
\label{lemma:newrolle} 
Suppose $C\!\subset\!\R^2$ is an arc (i.e., image of an interval or 
circle under a continuous map) with 
\begin{enumerate}
\item{A unique well-defined tangent line for each $x\!\in\!C$.}
\item{At most $I$ isolated\footnote{ Relative to the locus of 
inflection points.} inflection points.} 
\item{At most $V$ isolated points of vertical tangency.} 
\end{enumerate}
Then the maximum finite number of intersections of any line with $C$ 
is $I+V+2$. 
\end{lemma} 

\noindent 
{\bf Proof:} Let $S^1$ be the realization of the circle obtained by 
identifying $0$ and $\pi$ in the closed interval $[0,\pi]$. 
Consider the natural map $\phi : C\longrightarrow S^1$ obtained by 
$x\mapsto \theta_x$ where $\theta_x$ is the angle the normal line of 
$x$ forms with the $x_1$-axis. We claim that any $\theta\!\in\!S^1$ has at 
most $I+V+1$ pre-images under $\phi$. 

To see why, note that by assumption we can express $C$ as the union of 
no more than $I+V+1$ arcs where (a) any distinct pair of arcs is either 
disjoint or meets at $\leq\!2$ end-points, and (b) every end-point is either 
an isolated point of inflection or vertical tangency of $C$. Calling these arcs 
{\bf basic arcs}, it is then clear that the interior of any basic arc is 
homeomorphic (via $\phi$) to a connected subset of $S^1\!\setminus\!\{0\}$. 
Furthermore, by construction, the cardinality of $\phi^{-1}(0)$ is exactly 
$V$. So we indeed obtain that any 
$\theta\!\in\!S^1$ has at most $I+V+1$ pre-images under $\phi$. 

Now note that any line $\{x \; | \; m_1x_1+m_2x_2\!=\!m_0\}$ normal 
to $C$ forms an acute angle of  
$\mathrm{ArcTan(\frac{m_2}{m_1})}$ with the $x_1$-axis. Thus, the number of 
contact points $C$ has with the differential system 
\[ \frac{\partial x_1}{\partial t}=m_2 \ , \ \frac{\partial x_2}
{\partial t}=-m_1 \] 
is\footnote{i.e., the number of points at which some solution of 
the differential system has a tangent line in common with $C$ is...}
at most $I+V+1$. By Rolle's Theorem for Dynamical Systems in the Plane 
(see, e.g., \cite[corollary, pg.\ 23]{few}), we then obtain that 
the number of intersections of $\{x \; | \; m_1x_1+m_2x_2\!=\!m_0\}$ 
with $C$ is at most $I+V+2$, for any real $(m_0,m_1,m_2)$. So 
we are done. \qed 

\begin{rem} 
The bound from lemma \ref{lemma:newrolle} is tight in all cases. 
This is easily revealed by the following examples and their 
obvious extensions:\\ 
\vspace{-.5cm}
\begin{figure}[h]
\epsfig{file=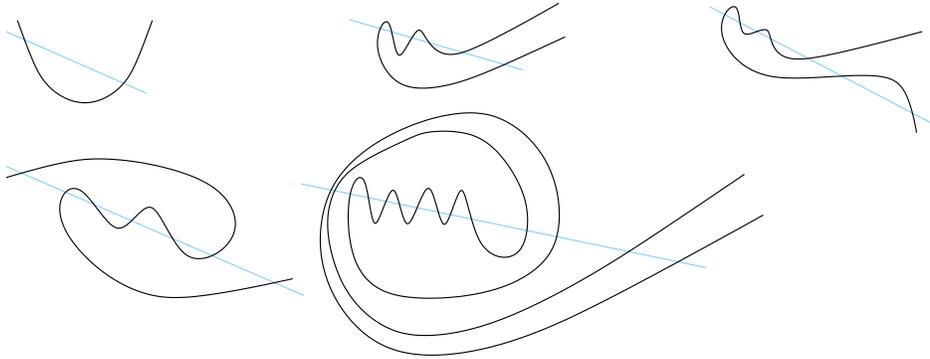,height=2in} 
\caption{Lemma \ref{lemma:newrolle} gives a tight bound for 
$(I,V)\!\in\!\{(0,0),(3,1),(4,1),(3,2),(7,5)\}$ and this generalizes easily to 
arbitrary $(I,V)$. } 
\end{figure} 

\vspace{-.5cm}
\noindent 
The authors do not presently know whether this bound remains tight when 
restricted to fewnomial zero sets. $\diamond$   
\end{rem}

We are now ready to give a quick geometrically motivated proof of 
the nearly optimal bound $\cN(3,3)\!\leq\!6$. 
This ``second'' proof of $\cN(3,3)\!\leq\!6$ was actually the original 
motivation behind this paper. 

\noindent
{\bf Short Geometric Proof of $\pmb{\cN(3,3)\!\leq\!6}$:}
Theorem \ref{thm:tri1} implies that we can assume that $f_1$ and $f_2$ have 
Newton polygons that are each triangles. Letting 
$Z$ denote the zero set of $f_2$ in $\R^2_+$, lemma \ref{lemma:tri2} of
the last section tells us that we can assume that $f_1\!=\!1\pm x_1\pm x_2$;
and by proposition \ref{prop:easy} the underlying change of variables also
implies that $Z$ is diffeomorphic to a line. So $Z$ is smooth and theorem 
\ref{thm:cool} tells us that $Z$ has no more than
$3$ inflection points and $1$ vertical tangent. So we now need only
check how many intersections $Z$ will have with the line $\{x \; | \;
1\pm x_1 \pm x_2\!=\!0\}$. By lemma \ref{lemma:newrolle}, we are done. \qed  

It turns out that inflection points for $m$-nomial curves are easy to describe 
in a $m$-nomial way. Let $\partial_i\!:=\!\frac{\partial}{\partial x_i}$. 
\begin{lemma} 
\label{lemma:imp} 
Suppose $f : \R^2_+ \longrightarrow \R$ is analytic and 
$Z$ is the real zero set of $f$. Then 
$z$ is an inflection point or a singular point of $Z \Longrightarrow 
f(z)\!=\!0$ and $[\partial^2_1 f\cdot  (\partial_2 f)^2- 
2\partial_1\partial_2 f\cdot\partial_1f\cdot\partial_2f 
+\partial^2_2 f\cdot(\partial_1 f)^2]_{x=z}\!=\!0$. 
In particular, in the case $f\!:=\!1+c_1x^{a_1}+\cdots +c_mx^{a_m}$,  
the preceding polynomial in derivatives is, up to a multiple which is a 
monomial in the $x^{a_i}$, a cubic polynomial homogeneous in the $c_ix^{a_i}$. 
\end{lemma}  
\noindent 
{\bf Proof:} In the case of a singular point, the first assertion 
is trivial. Assuming $\partial_2f\!\neq\!0$ at an inflection point then 
a straightforward computation of $\partial^2_1x_2$ (via implicit 
differentation and the chain rule) proves the first 
assertion. If $\partial_2f\!=\!0$ at an 
inflection point then we must have $\partial_1f\!\neq\!0$. So by computing 
$\partial^2_2x_1$ instead, we arrive at the remaining case of the first 
assertion. The second assertion also follows routinely. \qed 

Let us now reveal the hardest case of our result for pairs of trinomials. 
First note that while one can naturally associate a pair of polygons to 
$F$ when $n\!=\!2$, we can also associate a {\bf single} polygon by forming 
the {\bf Minkowski sum} $P_F\!:=\!\newt(f_1)+\newt(f_2)$. We can then give 
the following addendum to theorem \ref{thm:tri3} (with an independent proof).  
\begin{cor}
\label{cor:poly}  
Following the notation of GKC and theorem \ref{thm:tri3}, consider the case 
$(n,m_1,m_2)\!=\!(2,3,3)$.  Then $\cN(3,3)$ is respectively $0$,  
$2$, or $4$, according as we restrict to those $F$ with $P_F$ a line segment, 
triangle, or $\ell$-gon with $\ell\!\in\!\{4,5\}$. 
\end{cor} 

\noindent 
{\bf Proof:} The segment case follows immediately from 
corollary \ref{cor:zero}. For the remaining cases, proposition 
\ref{prop:easy} implies that we can assume $f_1\!:=\!1-x_1-x_2$ and 
$f_2\!:=\!1+Ax^a_1x^b_2+Bx^c_1x^d_2$. In particular, it is easily 
verified that the underlying monomial change of variables preserves 
the postivity of angles between lines (in exponent space), so the 
number of edges of $P_F$ is unchanged. 

Let $S_1\!:=\!Ax^a_1x^b_2$, $S_2\!:=\!Bx^c_1x^d_2$, and let 
$Z$ denote the zero set of $f_2$. Observe 
that lemma \ref{lemma:imp} tells us that we can bound the number 
of inflection points of $Z$ by analyzing the roots of a homogeneous 
polynomial in $(S_1,S_2)$ of degree $\leq\!3$. So let us now explicitly 
examine this polynomial in our polygonally defined cases. 

Clearly then, the triangle case corresponds to 
setting $a\!=\!d\!>\!0$ and $b\!=\!c\!=0$. 
We then obtain that [$x$ is an inflection point or a singular point of 
$Z$] $\Longrightarrow 1+S_1+S_2\!=\!0$ and $S_1+S_2\!=\!0$. So 
$Z$ has {\bf no} inflection points (or singularities). It is also even 
easier to see that $Z$ has no vertical tangents. So by lemma 
\ref{lemma:newrolle}, 
$\cN(3,3)\!\leq\!2$ in this case. To see that equality can hold in this 
case, simply consider $F\!:=\!(x^2_1+x^2_2-25,x_1+x_2-7)$, 
which has $P_F\!=\!\conv(\{(0,0),(3,0),(0,3)\})$ and root set 
$\{(3,4),(4,3)\}$.  

Similarly, the quadrilateral case corresponds to setting 
$b\!=\!c\!=\!0$ and $a,d\!>\!0$. We then get the pair of 
equations $1+S_1+S_2\!=\!0$ and $a(d-1)S_1-d(a-1)S_2\!=\!0$,  
with $a,d\!\not\in\{0,1\}$. (If $\{a,d\}\cap\{0,1\}\!\neq\!\emptyset$ 
then $F$, or a suitable pair of linear combination of $F$, would be  
pyramidal and we would be done by theorem \ref{thm:tri1}.)  
So $Z$ can have at most $1$ inflection point.  It is also even
easier to see that $Z$ has no vertical tangents. So by another application 
of lemma \ref{lemma:newrolle}, $\cN(3,3)\!\leq\!4$ in this case. To see that 
equality can hold in this case, simply consider the system 
$(x^2_1-3x_1+2,x^2_2-3x_2+2)$, which has 
$P_F\!=\!\conv(\{\bO,(2,0),(2,2),(0,2)\})$ and root set 
$\{(1,1),(1,2),(2,1),(2,2)\}$. 

Finally, the pentagonal case corresponds to setting $b\!=\!0$ and 
$a,c,d\!>\!0$. We then get the pair of equations $1+S_1+S_2\!=\!0$ and 
$a^2(d-1)S^2_1+a(ad-d-2c)S_1S_2-c(c+d)S^2_2\!=\!0$,   
with $ac(d-1)(c+d)\!\neq\!0$. (Similar to the last case, it is easily checked 
that if the last condition were violated, then we would be back in one of our 
earlier solved cases.) However, a simple check of the discriminant of the 
above quadratic form in $(S_1,S_2)$ shows that there is at most $1$ 
root, counting multiplicities, in any fixed quadrant. So, similar to the 
last case, we obtain $\cN(3,3)\!\leq\!4$ in this case. 
To see that the equality can hold in this case, simply 
consider the system $(x^2_2-7x_2+12,-1+x_1x_2-x^2_1)$, which has  
$P_F\!=\!\conv(\{\bO,(2,0),(2,2),(1,3),(0,2)\})$ and  
root set $\{(3,\frac{3\pm\sqrt{5}}{2}),(4,2\pm\sqrt{3})\}$. \qed  

\section{Monomial Morse Functions and Connected Components: Proving Theorem 
\ref{thm:cool}} 
\label{sec:morse}
A construction which will prove quite useful when we count connected 
components via critical points of maps is to find a monomial which is a  
Morse function relative to a given fewnomial zero set. 
\begin{rem} 
In what follows, we will always understand $\dim$ (resp.\ $\dim_\C$) to mean 
{\bf real} (resp.\ complex) dimension. Also, unless otherwise noted, 
``dimension'' will be understood to mean {\bf real} dimension. $\diamond$. 
\end{rem} 
\begin{lemma} 
\label{lemma:fiber}
Suppose $Z$ is the zero set in $\Rn_+$ of an $n$-variate $m$-nomial $f$. 
Then there exists a finite union of hyperplanes $H_Z\!\subset\!\Rn$ such that 
for all $a\!\in\!\Rn\!\setminus\!H_Z$ we have...
\begin{enumerate}
\item{Every critical point of the restriction of $x^a$ to $Z$ 
is non-degenerate.} 
\item{The level set in $Z$ of any regular value of $x^a$ has 
dimension $\leq\!n-2$.} 
\item{No connected component of $Z$ (other than an isolated point) is 
contained in any level set of $x^a$.}
\item{Every unbounded connected component of $Z$ has unbounded values 
of $x^a$.}  
\end{enumerate} 
\end{lemma}

\noindent 
{\bf Proof:} Let us prove the last two assertion first: Since the 
number of connected components of $Z$ is finite,\footnote{The smooth case 
is detailed in \cite[sec.\ 3.14]{few} and the case of integral exponents 
(allowing degeneracy) is a special case of \cite[lemma 3.2]{real}. In any 
event, the proof of the latter lemma extends easily to real exponents.} 
we can temporarily assume that $Z$ consists of a single connected component. 
Then, if we could find $n$ linearly independent $a$ with 
$Z\!\subset\!\{x\!\in\!\Rn_+ \; | \; x^a\!=\!c_a$\} for some $c_a$, 
proposition \ref{prop:easy} would immediately imply that $Z$ is contained 
in a point. Similarly, if we could find $n$ linearly independent $a$ 
for which the restriction of $x^a$ to $Z$ is bounded, then we would 
obtain by proposition \ref{prop:easy} again that $Z$ is bounded --- a 
contradiction. 

To prove the rest of our lemma, let us return to general $Z$ and consider the 
substitution $x_i\!=\!e^{z_i}$. A simple derivative computation (noting that 
$x\mapsto (e^{x_1},\ldots,e^{x_n})$ is a 
diffemorphism between $\Rn_+$ and $(\Rs)^n$) then shows that 
it suffices to instead prove the analogous statement where $f$ is replaced by 
a real exponential sum (a real analytic function in any event) and $x^a$ is 
replaced by the linear form $a_1z_1+\cdots+a_nz_n$. The latter analogue is then 
nothing more than an application of \cite[lemma 1, pg.\ 304]{bcss}, combined 
with Khovanski's Theorem on Fewnomials to ensure that $H_Z$ is finite instead 
of countable. \qed 

We will also need the following useful perturbation result, which 
can be derived via a simple homotopy argument. (See, e.g., 
\cite[lemma 2]{basu} for even stronger results of this form in the case of 
integral exponents.) 
\begin{lemma}
\label{lemma:pert} 
Following the notation of lemma \ref{lemma:fiber},  
let $Z_\delta$ denote the solution set of 
$|f|\!\leq\!\delta$ in $\R^n_+$ and 
$\stackrel{\circ}{Z_\delta}$ its boundary. Then for $\delta\!>\!0$ 
sufficiently small, $\stackrel{\circ}{Z_\delta}$ and its closure
are smooth, and there is a bijection between the connected components 
of $Z$ and $Z_\delta$ which preserves compact and non-compact components. \qed 
\end{lemma}

Finally, we will need the following two results 
(the latter dating back to an analogous result  
of Giusti and Heintz \cite[sec.\ 3.4.1]{giustiheintz} in the complex 
algebraic case, if not earlier) for dealing with over-determined 
fewnomial systems. 
\begin{dime} 
Suppose $U$ is an open subset of $\Rn$, $W$ is an irreducible real analytic 
subvariety  of $U$, and $g : U \longrightarrow \R$ is a real analytic function 
with $g(w)\!\neq\!0$ for some $w\!\in\!W$. Then 
$\dim W\cap \{z\!\in\!U \; | \; g(z)\!=\!0\}\!<\!\dim W$.  
\end{dime} 

\noindent 
{\bf Proof:} Let $d\!:=\!\dim W$ and let $W_\C$ be the complexification 
of $W$. Then $W_\C$ is an irreducible analytic subvariety  
of $U'$ where $U'\!\subseteq\!\Cn$ is an open subset containing $U$ and 
$\dim_\C W_\C\!\geq\!d$. Furthermore, by \cite[Active Lemma, pg.\ 
100]{gr} we have $\dim_\C W_\C\cap\{z\!\in\!U' \; | \; g(z)\!=\!0\}\!=\!
\dim_\C W_\C-1$. So, $W\cap \{z\!\in\!U \; | \;
g(z)\!=\!0\}$ (the real part of $W_\C\cap \{z\!\in\!U' \; | \; g(z)\!=\!0\}$)  
must have strictly smaller real dimension than $W$. \qed 
\begin{lemma}
\label{lemma:gh} 
Suppose $k\!\geq\!n$ and that $F\!:=\!(f_1,\ldots,f_k)$ is a $k\times n$ 
fewnomial system. Assume further that there are 
at most $m$ distinct exponent vectors in $F$. 
Then there exist real numbers $a_{ij}$ such that 
\begin{enumerate}
\item{the real zero 
set of $G\!:=\!(a_{11}f_1+\cdots+a_{1k}f_k,\ldots,a_{n1}f_1+\cdots+a_{nk}f_k)$ 
is the union of the real zero set of $(f_1,\ldots,f_k)$ and a finite 
(possibly empty) set of points.}
\item{$G$ is of type $(m-1,\ldots,m-1)$ and has no more than 
$m$ distinct exponent vectors. } 
\end{enumerate}
\end{lemma} 

\noindent
{\bf Proof:} Let us first make the substitution $x_i\!=\!e^{z_i}$, 
noting that $x\mapsto (e^{x_1},\ldots,e^{x_n})$ is a diffeomorphism between 
$\Rn$ and $\Rn_+$ which preserves the dimension of the underlying 
subanalytic varieties. Now pick $a_{11},\ldots,a_{1k}\!\in\!\R$ so that 
$a_{11}f_1+\cdots+a_{1k}f_k$ is not identically zero. 
Fix a set of points $\{w_i\}$, one lying in each 
irreducible component of the zero set $Z_1$ of $a_{11}f_1+\cdots+a_{1k}f_k$ 
in $\Rn$. Let us then pick $a_{21},\ldots,a_{2k}$ 
so that $a_{21}f_1+\cdots+a_{2k}f_k$ does not vanish at any 
$\{w_i\}$. By the Real Dimension Lemma we then obtain that the zero set $Z_2$ 
of $(a_{11}f_1+\cdots+a_{1k}f_k, a_{21}f_1+\cdots+a_{2k}f_k)$ in 
$\Rn$ is the union of a diffeomorphic copy of $Z$ and a real analytic variety 
of dimension $n-2$. Continuing this construction inductively, and then 
changing variables back again, we easily obtain assertion (1). 

An application of Gaussian elimination to eliminate one monomial 
from each of the polynomials of $G$ then gives us assertion (2). \qed 

To finally prove theorem \ref{thm:cool}, let us make one last 
definition. 
\begin{dfn} 
\label{dfn:big} 
Letting $f$ be a bivariate $m$-nomial and $Z$ the 
zero set of $f$ in the positive orthant, define...    
\begin{itemize} 
\item[$S(m):=$]{The maximal number of isolated singular points of such a $Z$.} 
\item[$I(m):=$]{The maximal number of isolated\footnote{Relative to the 
locus of inflection points.} inflection points 
of such a $Z$. } 
\item[$V(m):=$]{The maximal number of isolated\footnote{Relative to the 
locus of points of vertical tangency.} points of vertical tangency of $Z$. 
\qed } 
\end{itemize}
\end{dfn}

Theorem \ref{thm:cool} is then an immediate corollary of theorem 
\ref{thm:bounds} below. 
\begin{dfn} 
Let $\pmb{\cK(n,\mu)}$ be the maximal number of isolated roots in 
$\Rn_+$ of an $n\times n$ fewnomial system with exactly $\mu$ 
distinct exponent vectors. (So $\cK'(n,\mu)\!\leq\!\cK(n,\mu)$.) $\diamond$   
\end{dfn} 
\begin{thm} 
\label{thm:bounds} 
Theorem \ref{thm:cool} is true. Furthermore, defining 
$\cK(n,0)\!:=\!0$ and following the notation of definition 
\ref{dfn:big}, we also have the following inequalities: 
\begin{enumerate} 
\addtocounter{enumi}{3}
\item{$S(m),V(m)\!\leq\!\cK(n,m)$ } 
\item{$S(m)+I(m)\!\leq\!3\cK'(n,m)$ for $m\!\leq\!3$} 
\end{enumerate} 
\end{thm}

\noindent 
{\bf Proof:} Let us focus first on proving theorem \ref{thm:cool}: To prove 
assertions (1) and (2), note that we can 
divide by a suitable monomial so that $f$ has a nonzero constant 
term. By lemma \ref{lemma:pert}, we have that for $\delta\!>\!0$ sufficiently 
small, it suffices to bound the number of compact and non-compact connected 
components of $Z_\delta$ (a ``thickening'' of $Z$). In particular, 
$\stackrel{\circ}{Z}_\delta$, the boundary of $Z_\delta$, and its closure, 
can be assumed to be smooth. Noting that every connected component of 
$\stackrel{\circ}{Z}_\delta$ is contained in some connected component of 
$Z_\delta$, it then suffices to bound the number of connected components  
of $\stackrel{\circ}{Z}_\delta$. 

By proposition \ref{prop:easy} and lemma \ref{lemma:fiber}, we can 
pick an $n\times n$ matrix $A$ so that, after we make the change of variables 
$x\mapsto x^A$, the number of compact and non-compact real connected 
components of $\stackrel{\circ}{Z}_\delta$ is preserved and {\bf no} connected 
component of $\stackrel{\circ}{Z}_\delta$ of positive dimension is contained 
in a hyper-plane parallel to the $x_1$-coordinate hyperplane. Furthermore, we 
can also assume that every non-compact component of 
$\stackrel{\circ}{Z_\delta}$ has unbounded values of $x_1$. So we are 
now ready to use critical points to count connected components.  

Consider then the system of 
equations $G_\pm\!:=\!(f\pm\delta,x_2\partial_2f,\ldots,x_n\partial_nf)$, 
where $\partial_i$ denotes the operator $\frac{d}{dx_i}$. By construction, 
every compact connected component of $\stackrel{\circ}{Z_\delta}$ results 
in at least two extrema of the function $x_1$, i.e., 
$P_{\mathrm{comp}}(n,m)$ is bounded above by an integer no more than 
half of the total number of roots of $G_+$ and $G_-$.  (In particular, if $Z$ 
were smooth to begin with, then it would suffice to count the isolated roots 
of $G\!:=\!(f,x_2\partial_2f,\ldots,x_n\partial_nf)$ 
instead and omit the use of $Z_\delta$ and $G_\pm$.) Note also that by 
construction, all the roots of $G_\pm$ (or $G$) are non-degenerate. 
Furthermore, by a simple application of Gaussian Elimination, we obtain 
that $G_\pm$ (or $G$) is of 
type $(\underset{n}{\underbrace{m-1,\ldots,m-1}})$ 
(and there are no more than $m$ distinct monomial terms 
occuring in $G_\pm$ or $G$), so assertion (1) follows immediately. (The 
bound for $P_\mathrm{comp}(1,m)$ follows immediately from UGDRS.) 

To prove assertion (2) of theorem \ref{thm:cool}, another application of lemma 
\ref{lemma:fiber} (and our much used proposition \ref{prop:easy}) tells us 
that we can assume that every unbounded connected 
component of $Z$ has arbitrarily large values of $x_1$. 
For $\eps\!>\!0$ sufficiently small, we then observe that every such 
component induces at least one connected component of the intersection 
$Z'\!:=\!Z\cap\!\{x \; | \; x_1=\frac{1}{\eps}\}$. So fix an 
$\eps\!>\!0$ sufficiently small so that this holds for all 
unbounded components. (Recall that there are only finite many, cf.\ 
the proof of lemma \ref{lemma:fiber}.) Then, by substituting 
$x_1\!=\!\frac{1}{\eps}$ into $f$, we obtain a new fewnomial hypersurface 
$Z''\!\subseteq\!\R^{n-1}$, also defined by an $m$-nomial, with 
at least as many connected components as $Z$ has unbounded components. 
To conclude, note that under the change of variables $x\mapsto 
(x^{-1}_1,\ldots,x^{-1}_n)$, the bounded non-compact components of $Z$ are 
injectively embedded into the unbounded components of a new $m$-nomial 
hypersurface. So by what we've already proved for our unbounded components, 
we at last obtain $P_{\mathrm{non}}(n,m)\!\leq\!2(P_{\mathrm{comp}}(n-1,m)+
P_{\mathrm{non}}(n-1,m))$. 

The bound for $P_\mathrm{non}(2,m)$ then follows from the now classical 
moment map. That is, given any $n$-dimensional convex compact polytope 
$P\!\subset\!\Rn$, there is a real analytic diffeomorphism 
$\psi : \Rn_+ \longrightarrow \mathrm{Int}(P)$, where 
$\mathrm{Int}(P)$ denotes the interior of $P$ \cite[sec.\ 4.2]
{tfulton}. In particular, if one picks $P$ to be the Newton polygon 
of $f$ then there is a bijection between (a) the intersections of $\psi(Z)$ 
with the interior of an edge of $P$ with inner normal $w$, and (b) the roots 
of the initial term polynomial $\mathrm{in}_w(f)\!:=\!\sum c_ax^a$ in 
$(\Rs)\times\{1\}$, where the sum is over all $a\!\in\!\supp(f)$ 
with minimal inner product with $w$. Since any non-compact component $U$ of 
$Z$ results in $\psi(U)$ having at least $2$ intersections with the edges of 
$P$, UGDRS immediately implies our bound for $P_\mathrm{non}(2,m)$, 
not to mention our bound for $P_\mathrm{non}(1,m)$. (In fact, in our 
bound for $P_\mathrm{non}(2,m)$, we can even replace $m$ by the number of 
monomials corresponding to points on the boundary of $P$.) 

Assertion (3) of theorem \ref{thm:cool} follows immediately from 
assertion (4), which we will now prove. First note that the singular points of 
$Z$ are exactly the roots of the over-determined fewnomial system 
$F\!:=\!(f,x_1\partial_1f,x_2\partial_2f)$. By 
lemma \ref{lemma:gh} the singular points of $Z$ are also contained in the 
roots of the system $G\!:=\!(g_1,g_2)$, where $G$ is of type $(m-1,m-1)$, 
has no more than $m$ distinct exponent vectors, and 
each $g_i$ is a suitable linear combination of $f$, $x_1\partial_1f$, and 
$x_2\partial_2f$. Furthermore, the real zero set of $G$ is 
the union of the real zero set of $F$ and a (possibly empty) finite set of 
points. This proves the bound on $S(m)$, and the bound on $V(m)$ 
is proved in almost exactly the same way, starting with the polynomial 
system $(f,x_2\partial_2f)$ instead. So assertion (4) is proved. 

To prove assertion (5), note that by lemmata \ref{lemma:tri2} and 
\ref{lemma:imp}, $(x_1,x_2)$ is an inflection point 
or a singular point of $Z  
\Longrightarrow f\!=\!q\!=\!0$, where $q$ is a homogeneous polynomial, 
{\bf in the non-constant monomials terms of $\pmb{f}$}, of degree at most $3$. 
Letting $S_1,\ldots,S_{m-1}$ denote the non-constant monomials terms of $f$, 
note that each complex factor $q'\!:=\!\alpha_1S_1+\cdots+\alpha_{m-1}S_{m-1}$ 
of $q$ is a $j$-nomial for some $j\!\leq\!m-1$. (Note that 
the fundamental theorem of algebra tells us that $q$ indeed splits completely 
over $\C[S_1,\ldots,S_{m-1}]$, provided $m\!\leq\!3$.)  
Also note that if 
$\alpha_i\!\neq\!0$, the fewnomial systems $(1+S_1+\cdots+S_{m-1},q')$ and 
$G\!:=\!(1+S_1+\cdots+S_{m-1}-q'/\alpha_i,q')$ have 
the same zero set, and $\alpha_i$ must be nonzero for some $i$. However, $G$ 
is of type $(m-1,m-1)$, has no more than $m$ distinct exponent vectors, and 
has no degenerate roots. So the system $(f,q)$ has at most $3\cN'(m-1,m-1)$ 
isolated roots in the positive quadrant of the $(x_1,x_2)$-plane. So assertion 
(5) is proved. \qed 
\begin{rem} 
The equality $\cK'(n,m)\!=\!\cK(n,m)$ appears to be known only for 
$(n,m)\!\in\!(1\times \N)\cup\{(2,2),(2,3),(2,4)\}$ and $m\!=\!n+1$. These 
few cases follow easily from \mbox{theorems \ref{thm:tri3} and \ref{thm:tri1}, 
via remark \ref{rem:kho}. $\diamond$} 
\end{rem} 

\section*{Acknowledgements} 
The authors thank Alicia Dickenstein and Bernd Sturmfels for pointing 
out Haas' counter-example. Special thanks also go to Bertrand Haas 
for pointing out an error in an earlier version of lemma 
\ref{lemma:imp}, an anonymous referee for giving many nice corrections, 
and to Felipe Cucker, Jesus Deloera, Paulo Lima-Filho, and Steve Smale for 
some nice conversations. 

\footnotesize
\bibliographystyle{acm}

\end{document}